\documentclass[12pt]{amsart}
\usepackage{amssymb}

\usepackage{amsmath}
\usepackage{amsthm}
\usepackage{a4wide}
\usepackage{pdfsync}
\usepackage{hyperref}
\usepackage{enumerate}
\usepackage{color}
\usepackage{cancel}
\numberwithin{equation}{section}

\newtheorem{theorem}{Theorem}[section]

\newtheorem{lemma}[theorem]{Lemma}

\theoremstyle{remark}
\newtheorem{remark}[theorem]{Remark}

\theoremstyle{definition}

\def\eps{\varepsilon}

\def\rr{\mathbb{R}}

\def\th{\mathcal{T}_h}

\def\eps{\varepsilon}

\def\mh{\mathcal{H}}

\title[]{Sharp numerical approximation of the  Hardy  constant }

\author[L. I. Ignat]{Liviu I. Ignat}
\address{
    [1] Institute of Mathematics ``Simion Stoilow'' of the Romanian Academy, 21 Calea Grivitei Street, 010702 Bucharest, Romania.
    \newline\indent
[2] National University of Science and Technology Politehnica Bucharest, 313 Splaiul Independen\c tei, 060042 Bucharest, Romania.
    \newline\indent
[3] Academy of
Romanian Scientists, Ilfov Street, no. 3, Bucharest, Romania.}
\email[Corresponding author]{liviu.ignat@gmail.com}

\author[E. Zuazua]{Enrique Zuazua}
\address{
    [1] Friedrich-Alexander-Universit\"at Erlangen-N\"urnberg, 
    Department of Mathematics, Chair for Dynamics, Control, Machine Learning and Numerics (Alexander von Humboldt Professorship), 
    Cauerstr. 11, 91058 Erlangen, Germany.
    \newline\indent
    [2] Chair of Computational Mathematics,
 Deusto University, 48007 Bilbao, Basque Country, Spain.
    \newline\indent
    [3] Universidad Aut\'onoma de Madrid,
    Departamento de Matem\'aticas, 
    Ciudad Universitaria de Cantoblanco, 28049 Madrid, Spain.}
\email{enrique.zuazua@fau.de}
\date{}

\begin{document}

\begin{abstract}
We study the $P_1$ finite element approximation of the best constant in the classical Hardy inequality over bounded domains containing the origin in $\mathbb{R}^N$, for $N \geq 3$. 

Even though this constant is not attained in the natural Sobolev space $H^1$, our main result establishes an explicit, sharp, and dimension-independent rate of convergence proportional to $1/|\log h|^2$. 

The analysis combines an enhanced Hardy inequality with a remainder term involving logarithmic weights, refined approximation estimates for Hardy-type singular radial functions that serve as minimizing sequences, structural properties of continuous piecewise linear finite elements, and techniques from weighted Sobolev spaces.

We also consider other closely related spectral problems involving the Laplacian with singular quadratic potentials and obtain sharp convergence rates.
\end{abstract}

 
\maketitle

\section{Introduction and main results}

The Hardy inequality plays a central role in the analysis of partial differential equations (PDEs) with singular potentials, spectral theory, and the study of critical functional inequalities. In its classical form, it provides a lower bound on the Dirichlet energy in terms of a weighted $L^2$-norm involving the distance to the origin.

For bounded domains $\Omega \subset \mathbb{R}^N$ containing the origin and with $N \geq 3$, the best Hardy constant is known and sharp but notably not attained in the Sobolev space $H^1_0(\Omega)$. The same occurs in the whole space $ \mathbb{R}^N$. This lack of attainability poses significant challenges for the numerical approximation of the constant via variational methods.

In this paper, we study the convergence behavior of finite element approximations of the Hardy constant. Specifically, we focus on continuous piecewise linear  $P_1$ finite elements and establish an explicit convergence rate of order $1/|\log h|^2$ as the mesh size $h \to 0$. Our analysis relies on an improved Hardy inequality involving a remainder term with logarithmic weights, carefully constructed singular test functions, weighted Sobolev space estimates, and interpolation error bounds adapted to the singular nature of the Hardy inequality.

To be more precise, given a bounded domain $\Omega \subset \mathbb{R}^N$, $N \geq 3$,  containing the origin, let us consider the following minimization problem related to the Hardy inequality
\begin{equation}
\label{hardy.min}
 \Lambda_{N}(\Omega)=\inf_{u\in H^1_0(\Omega)} \frac{\|\nabla u\|^2_{L^2(\Omega)}}{\| \frac{u}{|x|}\|^2_{L^{2}(\Omega)}}.
\end{equation}
It is well known that this infimum is not attained, is  independent of the domain, and coincides with the Hardy constant in the whole space
\begin{equation}
\Lambda_{N}(\Omega)=\Lambda_N=\frac{(N-2)^2}4.
\end{equation}
The minimizer is not achieved in $H^1_0(\Omega)$. For instance, when $\Omega$ is a ball, the minimizer should be a radial function  of the form $u(r)=r^{-N/2+1}(a_1+a_2\log (r))$, which does not belong to $H^1_0(\Omega)$ but to a larger Hilbert space $\mathcal{H}$, which is essentially the closure of $H^1_0(\Omega)$ with respect to the norm  (see Section \ref{functional})
\begin{equation}
\| u\|^2_{\mathcal{H}}= \|\nabla u\|^2_{L^2(\Omega)}-\Lambda_{N}\left \|\frac{u}{|x|}\right\|_{L^2(\Omega)}^2.
\end{equation}

Given a
$P_1$ finite-element subspace $V_h$ of $H^1_0(\Omega)$ associated to a finite element mesh in $\Omega$, the Hardy constant can be approximated by the corresponding finite-dimensional minimization problem:\begin{equation}
\label{vh.constant}
   \Lambda_{h}(\Omega)=\inf_{u\in V_h} \frac{\|\nabla u\|^2_{L^2(\Omega)}}{\| \frac{ u}{|x|}\|^2_{L^{2}(\Omega)}}.
\end{equation}
The main result of this paper is the following theorem, which holds under standard assumptions on the finite element mesh, to be detailed below.
\begin{theorem}\label{main}
	Let  $\Omega$ be a  smooth, convex domain of $\rr^N$, $N\geq 3$,   and  $V_h$ be the space of $P_1$  finite elements  on $\Omega$. Then 
	\begin{equation}
\label{est.hardy.1}
  \Lambda_{h}(\Omega)-\Lambda_{N}\simeq \frac{1}{|\log h|^2}.
\end{equation}
\end{theorem}

This result fully clarifies an issue whose analysis was initiated in \cite{della2023finite}, where the one-dimensional case was treated.  In  \cite{della2023finite}, the rotational symmetry of the ball was exploited to derive convergence rates in specific configurations by reducing the analysis to the one-dimensional case.
The approach we introduce here is more flexible and applies to arbitrary smooth convex domains, using that the approximating minimizers have the same singularity at the origin as the radial ones. 

Note that the two-dimensional case is critical and an inverse square logarithmic correction of the Hardy inequality should be included \cite{MR1938711}. This would require further analysis.

Such results are now well established in other related contexts. In particular, the answer is well known for the classical Poincar\'e inequality, which is directly related to the first eigenvalue of the Dirichlet Laplacian. In that setting, the first continuous eigenvalue and its finite element approximation are known to be $h^2$-close \cite[Prop. 6.30, p. 315]{MR3405533}, \cite[Section 8, p. 700]{MR1115240}.

It is important to note, however, that the Poincar\'e inequality differs in two fundamental ways from the Hardy constant under consideration here. First, for the Poincar\'e constant, the infimum is actually attained, making it a true minimum. Second, this minimum is realized by the first eigenfunction of the Laplacian, which -- up to normalization -- is unique, belongs to $H_0^1(\Omega)$, and is in fact smooth.

In \cite{ignat2025}, the same issue was analyzed for the Sobolev constant, yielding polynomial convergence rates, with an order depending on both the dimension $N$ and the exponent $p$  in $W^{1, p}$. It is worth noting that the Hardy constant differs fundamentally from the Sobolev constant in that, while the Hardy constant is not attained, the Sobolev constant is achieved in the whole space, and the set of minimizers forms an explicitly characterizable, finite-dimensional (of dimension $N+2$) manifold.

Thus, the Hardy constant represents a new instance, with respect to the existing literature, in the sense that it is not achieved even in the whole space. This partially explains the logarithmic, rather than polynomial, approximation rate. The technical reason for the logarithmic rate becomes rather natural within the proof, which employs an improved Hardy inequality with a logarithmic correction  \cite[Th.~2.5.2]{MR2078192}, which replaces the Sobolev deficit estimates used in  \cite{ignat2025}.

We also consider some closely related spectral problems involving the Laplacian with a singular quadratic potential. 

Let us first consider
\begin{equation}
\label{hardy.min.2}
  \mu_1(\Omega)=\inf_{u\in H^1_0(\Omega)} \frac{\|\nabla u\|^2_{L^2(\Omega)}-\Lambda_{N} \|\frac{u}{|x|}\|_{L^2(\Omega)}^2}{\|  u \|_{L^{2}(\Omega)}}=\min_{u\in \mathcal{H}} \frac{\|\nabla u\|^2_{L^2(\Omega)}-\Lambda_{N} \|\frac{u}{|x|}\|_{L^2(\Omega)}^2}{\|  u \|_{L^{2}(\Omega)}}.
\end{equation}
Again, the optimal constant is not achieved in $H^1_0(\Omega)$ but in a larger Hilbert space $\mathcal{H}$. Of course, this spectral problem is well--posed thanks to the Hardy inequality, which ensures that the numerator of the Rayleigh quotient is non-negative and coercive in $\mathcal{H}$. The compactness of the embedding $\mathcal{H}\subset L^2(\Omega)$ assures that the minimum is achieved over $\mathcal{H}$.

We also consider its discrete counterpart, defining  $\mu_{1h}(\Omega)$ as the minimum of the same ratio in $V_h$.
In this case, we have the following result.
\begin{theorem}\label{main.2}In the setting of Theorem \ref{main}
	\begin{equation}
\label{est.hardy.2}
  \mu_{1h}(\Omega)-\mu_1(\Omega)\simeq \frac{1}{|\log h|}.
\end{equation}
\end{theorem}

Note that, although the convergence rate is also logarithmic, its order differs from that established in  Theorem~\ref{main}.

The same problem can be considered in the subcritical case in which the amplitude of the quadratic potential is $0\leq \Lambda<\Lambda_{N} 
$:
\begin{equation}
\label{hardy.min.3}
  \lambda_1(\Omega)=\min_{u\in H_0^1(\Omega)} \frac{\|\nabla u\|^2_{L^2(\Omega)}-\Lambda \|\frac{u}{|x|}\|_{L^2(\Omega)}^2}{\|  u \|_{L^{2}(\Omega)}}.
  \end{equation}
Note that, in this case, the minimizer is achieved in $H^1_0(\Omega)$ given the coercivity (in $H^1_0(\Omega)$) of the numerator of the Rayleigh quotient.

Minimizing over the finite element subspace $V_h$, we obtain the corresponding FEM approximation $ \lambda_{1h}(\Omega)$.
In the present case, we obtain a polynomial convergence rate:
\begin{theorem}\label{main.3}
In the setting of Theorem \ref{main}, denoting 
\begin{equation}m=\sqrt{\Lambda_N-\Lambda}\in \left(0,\frac{N-2}2\right],\end{equation}
 the following holds:

1. for $N\geq 5$
	\begin{equation}
\label{est.hardy.3}
  \lambda_{1h}(\Omega)-\lambda_1(\Omega)\simeq 
  \begin{cases}
h^{2m},& 0<m<1,\\
h^2|\log h|, & m=1,\\
h^2, &1<m\leq \frac{N-2}2;
  \end{cases}
\end{equation}

2. for $N=4$
\[
 \lambda_{1h}(\Omega)-\lambda_1(\Omega)\simeq 
\begin{cases}
h^{2m},& 0<m<1,\\
h^2, &m=1; 
  \end{cases}
\]

3. for $N=3$ 
\[
 \lambda_{1h}(\Omega)-\lambda_1(\Omega)\simeq 
\begin{cases}
h^{2m},& 0<m<\frac 12,\\
h^2, &m=\frac 12.
  \end{cases}
\]
\end{theorem}

In the present setting, the continuous eigenfunctions are less singular than in the critical case $\Lambda=\Lambda_N$, and belong to  $H_0^1(\Omega)$, although they still exhibit a singularity of order  $|x|^{-N/2+1+\sqrt{\Lambda_N-\Lambda}}$. Thus, the error analysis reduces to estimating the $H_0^1(\Omega)$--distance between the eigenfunction and the finite element space.  This permits obtaining polynomial decay rates rather than logarithmic ones, as in the previous theorems. 

Note that in the case $\Lambda=0$, we recover the classical problem of finite element approximation of the Poincar\'e constant, and we retrieve the well--known optimal convergence rate of order $h^2$.

The paper is organized as follows. We first introduce the functional framework needed to address the Hardy inequality, along with some classical results of finite element theory. We then present the proofs of the main results. The paper concludes with a section on comments and open problems, followed by an appendix containing technical lemmas.

\section{Functional framework}\label{functional}

As mentioned above, we denote by $\mathcal{H}\subset L^2(\Omega)$, the Hilbert space obtained as the completion of $C_c^\infty(\Omega)$ with respect to the norm
\begin{equation}
\label{hardy.norm}
  \|u\|_{\mathcal{H}}^2=\int_{\Omega} \left(|\nabla u|^2-\Lambda_{N}\frac{u^2}{|x|^2}\right)dx.
\end{equation}
The space $\mathcal{H}$ is isometric to the space 
\[
W^{1,2}_{0}(|x|^{-(N-2)}dx,\Omega)=W^{1,2}(d\mu, \Omega), \|v\|_{\tilde{H}}=\left(\int_{\Omega} |x|^{-(N-2)}|\nabla v|^2dx\right)^{1/2} 
\]
under the transformation 
\[u=Tv=|x|^{-( N-2)/2}v.\]
Here and in what follows, the density $\mu$ is defined as
\begin{equation}
\mu(x)= |x|^{-(N-2)}.
\end{equation}
It has been proved in \cite{MR1760280,MR2960851} that 
\[
\mathcal{H}=\overline{\{u=Tv, \ v\in C_c^\infty(\Omega)\}}^{\|\cdot\|_{\mathcal{H}}},
\]
where 
\[
\|u\|_{\mathcal{H}}^2=\int_{\Omega} |x|^{-(N-2)}|\nabla (|x|^{(N-2)/2}u)|^2dx.
\]

Let us introduce the two operators 
\[
Lu=-\Delta u -\Lambda_{N} \frac{u}{|x|^2},
\]
\[
\tilde Lv=-|x|^{N-2}\nabla \cdot (|x|^{-(N-2)}\nabla v)=-\Delta v+(N-2)\frac{x}{|x|^2}\cdot \nabla v.
\]
Under the transformation $u=Tv$, we have
\[
Lu=L(Tv)=|x|^{-N/2+1}\tilde Lv=T\tilde Lv, 
\]
i.e.  $LT=T\tilde L$. 

There exists a sequence of common eigenvalues 
\[
0<\mu_1\leq \mu_2\leq \cdots \leq  \rightarrow\infty
\]
and the corresponding eigenfunctions are related by
$\phi_k=T\psi _k$. 

In the particular case when $\Omega$ is a ball centered at the origin (for simplicity, we take $\Omega=B_1(0)$), the eigenvalues can be computed explicitly 
\[
\psi_{j,n}(r,\sigma)=J_{m_j}(z_{m_j,n}r)f_j(\sigma),
\]
where $\{f_j\}_{j\geq 0}$ are the spherical harmonics, $z_{m,n}$ is the $n$th zero of the Bessel function $J_m$, $m^2=j(j+N-2)$, $j\geq 0$, and $\mu_{j,n}=z_{m_j,n}^2$.

Let us denote by $\psi_1$ the first eigenfunction, i.e., 
\[
\psi_1=\psi_{0,1}=J_0(z_{0,1}r),
\] 
and $\phi_1=T \psi_1$. The properties of the Bessel function guarantee that 
\[
\psi_1(0)\neq 0, \ \nabla \psi_1(0)=0.
\]

\section{Fundamental Tools from Finite Element Analysis}
\label{finite.elements}
\subsection{Finite element meshes}
Let  $\Omega$  be a polytope in $\rr^N$  (in particular, an interval, a polygonal or a polyhedral domain in dimension one, two, and three, respectively). For each positive $h$ we construct a partition $\th$ (a mesh) of the domain $\Omega$ into a finite set of $N$-simplices or cells (tetrahedra in  dimension $N=3$) satisfying 
\begin{enumerate}
\item $\cup _{T\in \th} T=\overline{\Omega}$,
\item if $T,T'\in \th$, $T\neq T'$, then either $T\cap T^{'}=\emptyset $, $T\cap T^{'}$ is a single common vertex or $T\cap T^{'}$ is a whole common facet (point in $N=1$, edge in $N=2$, face in $N=3$). 
\end{enumerate}

For each $T\in \th$, we denote by $\rho_T$ and $h_T$  the radius of the largest inscribed ball in $T$ and the diameter of $T$, respectively. We set 
\[
h=h(\th)=\max_{T\in \th}h_T.
\] 
We will consider a set of regular meshes $(\th)_{h>0}$:  there exists $\sigma>0$, independent of $h$, such that
\begin{equation}
\label{regular}
  \frac{h_T}{\rho_T}\leq\sigma, \ \forall \ T\in \th, \ \forall\ h>0.
\end{equation}
The mesh is also assumed to be quasi-uniform, i.e. 
\[
\inf_{h>0}\frac {\min_{T\in \th}h_T}{ \max _{T\in \th} h_T}>0.
\]
 
Each element of the mesh $\th$ is the image of a reference $N$-simplex  through an affine mapping $F_T:\rr^N\rightarrow \rr^N$, 
\[
F_T(\widehat x)=B_T \widehat x+b_T,
\]
with $B_T$ being an invertible $N\times N$ matrix, $b_T\in \rr^N$, such that
\[
F_T(\widehat T)=T, \quad \ \forall \ T\in \th.
\]
We recall a few properties of matrix $B_T=\nabla F_T$ (cf. \cite[Lemma 7.4.3, p.~272]{buttazzo})
\[
\|B_T\|\leq \frac{h_T}{\rho_{\widehat T}}, \quad \|B_T^{-1}\|\leq \frac{h_{\widehat T}}{\rho_{T}}
\]
and 
\[
|\det B|=|  J F_T|=\frac{|T|}{|\widehat{T}|}.
\]

For a fixed $\th$, we define  the space $\widetilde V_h$ as 
\[
\widetilde V_h=\{ f\in C(\overline\Omega)\cap H_0^1(\Omega); f\circ F_T\in \mathbb{P}^1(\widehat T), \ \forall T\in \th \},
\]
where $ \mathbb{P}^1(\widehat T)$ is the space of linear polynomials on $\widehat T$.

For our approximations, we consider   the unit ball $
 B$, and we approximate it with a polygonal domain $B_h\subset B$ as in \cite{MR2754218}.  We also introduce the space $V_h\subset H_0^1(B )\cap C(\overline B )$ of functions in $\widetilde V_h$ extended by zero in $B\setminus B_h$.

\subsection{Approximation by piecewise linear functions}

We recall the classical approximation result in Sobolev spaces \cite[Th. 4.4.20]{MR2373954}: For any polyhedral domain $\Omega\in \rr^N$,   $N<2p$, $\mathcal{T}_h$ as above, $1<p\leq \infty$ or $N\leq 2$ if $p=1$ (see  \cite[Th. 4.4.4]{MR2373954} for the complete set of restrictions), and $s\in \{0,1\}$, the global piecewise-linear interpolant $I^h$ satisfies  
 \begin{equation}
\label{fem.error.sobolev}
  \left(\sum_{T\in \mathcal{T}_h }\| u-I^h u\|^p_{W^{s,p}(T)} \right)^{1/p}\leq Ch^{2-s}\|u\|_{W^{2,p}(\Omega)}, \ \forall u\in W^{2,p}(\Omega).
\end{equation}
  The above restrictions are necessary when one needs to have an estimate for all functions in $W^{2,p}(\Omega)$.   When we restrict to the class of functions that are smooth, for example $C^2(\overline{\Omega})$, the above restrictions on the dimension can be relaxed. The restriction $N<2p$ is exactly the one that guarantees that $W^{2,p}(\rr^N)\subset C^0(\Omega)$, the class where the global linear interpolator $I^h$ is defined. 

\begin{lemma}\label{est.c2}Let $\mathcal{T}_h$ be a regular mesh on a 
	polyhedral domain $\Omega\in \rr^N$,    $1<p\leq \infty$, $s\in \{0,1\}$. There exists a positive constant $C=C(N,p,s,\sigma)$ such that for all $|\alpha|$=s and $1\leq p<\infty$:
	\begin{equation}
\label{est.fem.c2}
    \left(\sum_{T\in \mathcal{T}_h }\| D^\alpha(u-I^h u)\|^p_{L^{p}(T)} \right)^{1/p}\leq Ch^{2-s}\left( \sum_{T\in \mathcal{T}_h} |T| \|D^2 u\|^p_{L^\infty(T)}\right)^{1/p}, \ \forall u\in C^{2}(\overline{\Omega}), 
\end{equation}
and
	\begin{equation}
\label{est.fem.c2.infty}
    \max_{T\in \mathcal{T}_h}\| D^\alpha(u-I^h u)\|_{L^{\infty}(T)}\leq Ch^{2-s}\max_{T\in \mathcal{T}_h}   \|D^2 u\|_{L^\infty(T)}, \ \forall u\in C^{2}(\overline{\Omega}).
\end{equation}
\end{lemma}

The proof is a slight modification of the one in \cite[Th.~4.4.4, Chapter 4]{MR2373954}. The proof can be extended to more general finite elements and $u\in C^m(\overline{\Omega})$, but this is beyond the scope of this paper. For a proof, see \cite{ignat2025}.

In fact, the above estimates are sharp in the case of functions that are uniformly convex in one direction.  

\begin{lemma}\label{lower.bound.triangle}{\rm{(\cite{ignat2025})}}For any $p\in (1,\infty)$ there exists a positive constant $C(p)$ such that for any $T\in \mathcal{T}_h$ and any
 $u\in C^2(T)$  
	\[
	 \min_{A\in \rr^N} \int_T |Du-A|^pdx \geq C(p)\rho_T^{N+p} \max _{\xi\in \mathbb{S}^{N-1}}\min _{x\in \overline{\Omega}}|\xi^T D^2u(x)\xi|^p.
	\]
	The same holds under the assumption that the function is uniformly convex in one direction, i.e. $\inf_{T}|\partial^2_{x_k}u(x)|>0$ for some $x_k$.
	\end{lemma}

\subsection{Finite element eigenvalue approximation}\label{classical.approximation} Let us now recall the classical theory for eigenvalue approximation, \cite{MR1115240}. Here we present it in the simplest case.  Following the notations in \cite[Section 8, p.~697]{MR1115240}, let $V$ be a real Hilbert space and $a(\cdot, \cdot)$ be a symmetric continuous and coercive bilinear form on $V $. Let $H$ be another Hilbert space such that $V\subset H$ with compact embedding, $b$ a symmetric continuous bilinear form on $H$, such that $b(u,u)>0$, for all $u\in V$, $u\neq 0$. Let $V_h\subset V$ be a family of finite-dimensional spaces of $V$. 

Let $\lambda_1$ be the first eigenvalue of \textit{the form $a$ relative to the form $b$}, i.e. the smallest $\lambda_1$ so that there exists a non--trivial $u_1\in  V$ such that
\[
a(u_1,v)=\lambda_1b(u,v), \forall \ v\in V.
\]
In a similar way, we define $\lambda_{1h}$ as the smallest value for which there exists a non--trivial $u_{1h}\in V_h$ such that 
\[
a(u_{1h},v_h)=\lambda_{1h}b(u_{1h},v_h), \forall \ v_h\in V_h.
\]
A fundamental result in the theory of eigenvalue approximation is the following, originally established in \cite[Prop.~6.30, p.~315]{MR3405533}; see also \cite[p.~700]{MR1115240} for a more direct statement:
\begin{equation}
\label{two-side.estimate}
  C_1\eps_h^2 \leq \lambda_{1h}-\lambda_1\leq C_2\eps_h^2
\end{equation}
where
\[
\eps_h=d(u_1,V_h)=\inf_{v_h\in V_h } \|u_1-v_h\|_{V}. 
\]
The upper bound together with the trivial estimate $0\leq \lambda_{1h}-\lambda_1$,  can also be found in \cite[Th. 6.4-2]{raviart}.

This analysis will be useful in the proofs of Theorem~\ref{main.2} and Theorem~\ref{main.3}, which are effectively related to eigenvalue problems. However, Theorem~\ref{main} lies outside this framework and requires an independent, substantially more advanced analysis.

\section{Proof of Theorem \ref{main} }\label{hardy}
Let us consider $\Omega$ a convex smooth domain such that 
$0\in \Omega \subset \overline{\Omega}\subset B_R$ for some $R>0$.  Without loss of generality, assume that $R=1$. 
We approximate the domain $\Omega$ by a polygonal domain $\Omega_h\subset \Omega$ as in \cite{MR2754218}. We introduce the space $V_h$ to be the space of continuous finite elements on $\Omega_h$ extended by zero 
   to $\Omega\setminus \Omega_h$. 

The proof of Theorem \ref{main} treats the lower and upper bounds in \eqref{est.hardy.1} separately. 

The following Hardy inequality with a logarithmic remainder term \cite[Th.~2.5.2, p.~25]{peral} will play a key role: There exists a positive constant $C=C(N,\Omega)$ such that for all $\phi\in C_c^\infty(\Omega)$, it holds
\begin{equation}
\label{hardy.log}
  \int_\Omega |\nabla \phi|^2 dx- \Lambda_N\int_{\Omega}\frac{\phi^2}{|x|^2}dx\geq C \int_{\Omega} |\nabla \phi|^2 \left(\log\frac 1 {|x|} \right)^{-2}dx.
\end{equation}
By density, it also holds for any $\phi\in \mh$. We also recall that   (see \cite[(2.5.7), p.~25-26]{peral}) 
\begin{equation}
	\label{peral1.4}
\inf_{\phi\in C_c^\infty(\Omega), \phi\not= 0} \frac{  \int_\Omega |\nabla \phi|^2 dx- \Lambda_N\int_{\Omega}\frac{\phi^2}{|x|^2}dx}{ \int_{\Omega} | \phi|^2 |x|^{-2} \left(\log\frac 1 {|x|} \right)^{-2}dx}=\frac 14.
\end{equation}
In dimension $N\geq 3$, the infimum above, $1/4$, can be approximated through the  procedure in  \cite[Chapter~2.4, p.~22]{peral} and the fact that
\begin{equation}
\label{tildee2}
  \tilde u_2(x)=|x|^{-\frac{N-2}2} \left(\log \frac 1{|x|}\right)^{1/2} 
\end{equation}
is the distributional solution of 
\[
-\Delta w - \Lambda_N \frac w{|x|^2}=\frac 14 \frac {w}{|x|^2} \left(\log \frac 1{|x|}\right)^{-2}.
\]
This regularization consists, roughly speaking, of truncating the singular explicit function  $\tilde u_2$ 
near the singularity at $x=0$.

  It is worth mentioning that functions of the form $u(x)=|x|^{-(N-2)/2}v(x)$, where $v(x)\simeq (\log (1/x))^\alpha$ as $x\rightarrow  0$ belong to the space ${\mathcal{H}}$ if $\alpha\in [0,1/2)$ but fail to be in ${\mathcal{H}}$ for $\alpha\geq 1/2$, see \cite{MR2960851}.

The idea of the proof is to construct      a minimizing sequence $u_\eps$ for  $\Lambda_N$, with a controlled error of order $|\log \eps|^{-2}$: 
\[
\frac{  \int_\Omega |\nabla u_\eps|^2 dx- \Lambda_N\int_{\Omega}{u_\eps^2}{|x|^{-2}}dx}{ \int_{\Omega} | u_\eps|^2 |x|^{-2} }
\lesssim |\log \eps|^{-2}.
\]
  The proof of Theorem \ref{main} proceeds by a careful finite element approximation of this minimizing sequence.

\subsection{The lower bound.} Although it is briefly presented in \cite{della2023finite}, we include a sketch here for completeness.

 Let $v_h\in  V_h$ be the minimizer corresponding to $\Lambda_h$:
\[
\Lambda_h=\frac{\int _\Omega |\nabla v_h|^2dx}{\int _\Omega v_h^2 |x|^{-2}dx}. 
\]
It follows that
\begin{align*}
\Lambda_h-\Lambda_N= \frac{\int _\Omega |\nabla v_h|^2dx- \Lambda_N\int _\Omega v_h^2 |x|^{-2}dx }{\int _\Omega v_h^2 |x|^{-2}dx}\geq C(N,\Omega)  \frac{\int _\Omega |\nabla v_h|^2 |\log |x||^{-2}dx}{\int _\Omega v_h^2 |x|^{-2}dx}.
\end{align*}
Using that $\nabla v_h$ is constant in each simplex and  Lemma \ref{logth} in the Appendix below, we get
\begin{align*}
\label{}
  \int _\Omega |\nabla v_h|^2 |\log |x||^{-2}dx&=\sum_{T\in \th} \int_{T} |\nabla v_h|^2 |\log |x||^{-2}dx=\sum_{T\in \th} |\nabla v_h|^2\int_{T}  |\log |x||^{-2}dx\\
  &\gtrsim \frac{1}{|\log h|^2}\sum_{T\in \th} |T| |\nabla v_h|^2  =
 \frac{1}{|\log h|^2}\int_{\Omega}|\nabla v_h|^2dx\\
 &=\frac{\Lambda_{h}}{|\log h|^2}\int _\Omega v_h^2 |x|^{-2}dx.
\end{align*}
This shows that $\Lambda_h-\Lambda_N\gtrsim\Lambda_h|\log h|^{-2}\geq \Lambda_N|\log h|^{-2}$.  

\subsection{The upper bound}\label{upperbound}  It is sufficient to prove the upper bound 
\[
\Lambda_{h} -\Lambda_N\lesssim \frac 1{|\log h|^2},
\]
in the case when $\Omega=B_1$,  the unit ball. Indeed, given that the constant $\Lambda_N$ is independent of the bounded domain under consideration, the upper bound holds in the general case by comparison.


As mentioned above, the main idea for the upper bound is to use an approximating sequence $u_\eps$ of 
 \begin{equation}\label{u2alpha}
  \tilde u_2(x)=|x|^{-N/2+1}\left(\log \frac{1}{ |x|}\right)^\alpha,
\end{equation} 
 and project it on the space $ V_h$.  In the following, we make this construction more precise. 
 
In the proof, we will also employ the function
 \[  u_2(x)=|x|^{-N/2+1}.\]

We consider the following cut-off function inspired by \cite[proof of Corol. VIII.6.4]{MR3823299} and \cite{MR3528527}:
\[
\eta_\eps (x)=\begin{cases}
	0,& |x|<\eps^2, \\
	\xi (\frac {\log (|x|/\eps^2)}{\log (1/\eps)}), & |x|\in (\eps^2,\eps),\\
	1, & |x|>\eps,
\end{cases}
\]
where $\xi:[0,1]\rightarrow [0,1]$ is a smooth function such that for some $\mu\in (0,1)$,
 $\xi=0$ on $[0,\mu]$ and $\xi=1$ on $[1-\mu,1]$. 

The function $\eta_\eps$ satisfies the following properties:
 \begin{enumerate}
 \item \label{p0} $\eta_\eps(r)=0$ on $0<r<\eps^{2-\mu}$,  
 	\item \label{p1} $|\eta'_\eps(r)|\leq \frac{\|\xi \|_\infty}{r |\log \eps|}\lesssim \frac{1}{r|\log \eps|}$, $\eps^2\leq r\leq \eps$, 
 	\item \label{p2} $|\eta''(r)|\leq \frac{\|\xi''\|_\infty}{r^2|\log \eps|^2}+\frac{\|\xi'\|_\infty}{r^2|\log \eps|}\lesssim \frac{1}{r^2 |\log \eps|}$, $\eps^2<r<\eps$.
 \end{enumerate}
 As a consequence, we get:
 \begin{enumerate}
 	\item \label{p3}  $|\nabla \eta_\eps(x)|\leq \frac 1{|x||\log \eps|}$ if $\eps^2<|x|<\eps$ and vanishes otherwise,
 	\item \label{p5} $|D^2\eta_\eps (x)|\lesssim \frac 1{|x|^2 |\log \eps|}$, if $\eps^2<|x|<\eps$ and vanishes otherwise.
 \end{enumerate}
With this function $\eta_\eps$, we introduce 
\[
u_\eps(x)= \tilde u_2(|x|) \eta_\eps(|x|)\psi(|x|)\in C_c^\infty(\Omega),
\]
where $\psi\in C_c^\infty (\rr)$ such that $\psi \equiv 1$ for $|r|\leq 1/4$ and $\psi\equiv 0$ for $|r|>1/2$.

The following lemma provides quantitative estimates for this sequence, viewed as a minimizing sequence for the Hardy constant.

\begin{lemma}\label{est.u.eps}Let $\alpha\geq 0$. The family of functions $u_\eps$ satisfies the following estimates:
 \begin{equation}
\label{estAeps}
  A_\eps= \int_\Omega |\nabla u_\eps|^2dx-\Lambda_N\int_\Omega \frac{u_\eps^2}{|x|^2}dx\lesssim 
  \begin{cases}
  |\log \eps|^{2\alpha-1}, &\alpha > 1/2,\\[10pt]
\log |\log \eps|,& \alpha=1/2,\\[10pt]
1,& \alpha\in [0,1/2); 
 \end{cases}
 \end{equation}
	 \begin{equation}
\label{estBeps}
  B_\eps=\int _{\Omega} \frac{u_\eps^2}{|x|^2}dx\simeq |\log \eps|^{2\alpha+1},
  \end{equation}
 and  
\begin{equation}
\label{est.H2u.eps}
  \|u_\eps\|^2_{H^2(\Omega)}\lesssim   \frac{|\log \eps|^{2\alpha}}{\eps^{4-2\mu}}.
\end{equation}
\end{lemma}

\begin{remark}
	As a consequence, we get that the quotient $A_\eps/B_\eps$ satisfies
	\begin{equation}
\label{a/b}
\frac {A_\eps}{B_\eps}\lesssim 
\begin{cases}
	\frac 1{|\log \eps|^2},& \alpha>\frac 12,\\[10pt]
	\frac{\log |\log \eps|}{|\log \eps|^2},& \alpha=\frac 12,\\[10pt]
	\frac{1}{|\log \eps|^{2\alpha+1}}, & \alpha\in [0,\frac 12).
\end{cases}  
\end{equation}
\end{remark}

 \begin{proof}  Using that $u_\eps(x)=\tilde u_2(|x|) \eta_\eps(r)\psi(r)$,  
 we obtain
\[
B_\eps\leq \int_{\eps^2}^{1/2} r^{N-3}|\tilde u_2(r)|^2 \eta^2_\eps(r)dr\lesssim \int_{\eps^2}^{1/2} \frac{|\log r|^{2\alpha}}r  dr\lesssim |\log \eps|^{2\alpha+1}
\]
and
\[
B_\eps\gtrsim \int_{\eps}^{1/4} \frac{|\log r|^{2\alpha}}r dr\gtrsim |\log \eps|^{2\alpha+1},
\]
and then \eqref{estBeps} holds. 

The estimate \eqref{estAeps} is more delicate. 

Let us recall that for any $u\in \mathcal{H}$, 
\[
\|u\|_{\mh}^2=\int_{\Omega} |x|^{-(N-2)} \left|\nabla \left(\frac {u}{u_2}\right)\right|^2dx=\int_{\Omega}|\nabla u|^2dx-\int_{\Omega}\Lambda_N\frac{u^2}{|x|^2}dx.
\]
For
$u_\eps$, we have  \[
A_\eps =\int _{\Omega} \left| \nabla u_\eps-u_\eps \frac{\nabla u_2}{u_2} \right|^2dx=  \int_\Omega |u_2\nabla \theta_\eps|^2dx,
\]
where \begin{equation}\theta_\eps(r)=\left(\log \frac 1r\right)^{\alpha}\eta_\eps(r)\psi(r).\end{equation}

Note that 
\[\nabla \theta_\eps=\left (\left(\log \frac 1{|x|}\right)^{\alpha}\nabla \eta_\eps - \alpha\left(\log \frac 1 {|x|}\right)^{\alpha-1}\frac {x}{2|x|^2}\eta_\eps\right)\psi+ \left( \log \frac 1{|x|}\right)^\alpha \eta_\eps\nabla \psi.
\]
 According to the properties of $\psi$ and $\eta_\eps$, 
\begin{align*}
   \int_{\Omega}|u_2\nabla \theta_\eps|^2dx& \lesssim  \int_{|x|<1/2}\left(u_2^2 |\nabla \eta_\eps|^2 |\log |x||^{2\alpha}+u_2^2\eta_\eps^2 \frac{|\log |x||^{2\alpha-2}} {|x|^2 }\right)dx+\int _{1/4<|x|<1/2}|\log |x||^{2\alpha}\eta_\eps^2dx \\
   &\lesssim\int_{\eps^2}^{\eps } r|\log r|^{2\alpha}\frac{dr}{r^2 |\log \eps|^2} + \int_{\eps^2}^{1/2} \frac{|\log r|^{2\alpha-2}dr}{r} +\int_{1/4}^{1/2}  r{|\log r|^{2\alpha}dr}\\
   &\lesssim 1+\begin{cases}
|\log \eps|^{2\alpha-1}, &\alpha > 1/2,\\
 \log |\log \eps|,& \alpha=1/2,\\
 |\log \eps|^{2\alpha-1},& \alpha\in [0,1/2).
 \end{cases}
\end{align*}
This completes the proof of \eqref{estAeps}.

Let us now estimate the $H^2(\Omega)$ norm of $u_\eps$. For simplicity, assume that $\psi \equiv 1$, since the main contribution comes from the singularity near the origin. As above, we have
\[
\|u_\eps\|^2_{L^2(\Omega)}=\int_{\eps^2}^{1/2} r^{N-1}|\tilde u_2(r)|^2 \eta^2_\eps(r)dr\lesssim \int_{\eps^2}^{1/2} r {|\log r|^{2\alpha}}  dr\lesssim 1.
\]
Using the expression of $u_\eps$ and the properties of $\eta_\eps$ we have
\begin{align*}
\label{}
  |u_\eps'(r)|&\lesssim r^{-N/2}|\log r|^{\alpha}\eta_\eps(r)+r^{-N/2}|\log r|^{\alpha-1}\eta_\eps(r)+r^{-N/2+1}|\log r|^{\alpha} \eta'_\eps(r)\\
  &\lesssim r^{-N/2}|\log r|^{\alpha}\eta_\eps(r)+r^{-N/2+1}|\log r|^{\alpha} \eta'_\eps(r)\\
  &\lesssim r^{-N/2}|\log r|^{\alpha}+r^{-N/2}\frac{|\log r|^{\alpha}}{|\log \eps|}\lesssim r^{-N/2}|\log r|^{\alpha}.
\end{align*}
This implies that
\begin{align*}
  \|\nabla u_\eps\|_{L^2(\Omega)}^2&=\int_{\eps^{2-\mu}}^{1/2} r^{N-1}|u_\eps'(r)|^2dr\lesssim 
\int_{\eps^{2-\mu}}^{1/2} r^{-1}|\log r|^{2\alpha}dr\lesssim |\log \eps|^{2\alpha+1}.
\end{align*}
The second order derivatives of $u_\eps $ satisfy
\begin{align*}
    |D^2u^\eps(x)| \sim |u_\eps''(r)| \lesssim  & \quad  \eta_\eps(r)r^{-N/2-1}\left(|\log r|^{\alpha}+|\log r|^{\alpha-1}+|\log r|^{\alpha-2} \right)  \\
  &+\eta_\eps'(r)r^{-N/2}\left(|\log r|^{\alpha}+|\log r|^{\alpha-1} \right)\\
  &+\eta_\eps''(r)r^{-N/2+1}|\log r|^{\alpha}\\
  \lesssim  &\quad  \eta_\eps(r)r^{-N/2-1}|\log r|^{\alpha}+\eta_\eps'(r)r^{-N/2} |\log r|^{\alpha}+\eta_\eps''(r)r^{-N/2+1}|\log r|^{\alpha}.
\end{align*}
This implies that
\begin{align*}
\int_\Omega |\partial_{ij}u_\eps|^2dx&\lesssim \int_{\eps^{2-\mu}}^{1/2}r^{N-1}\left(|u_\eps''(r)|^2+\frac{|u_\eps'(r)|^2}{r^2}\right)dr\\
&\lesssim \int_{\eps^{2-\mu}}^{1/2} \left(|\eta_\eps(r)|^2r^{-3}|\log r|^{2\alpha}+|\eta_\eps'(r)|^2r^{-1}|\log r|^{2\alpha}+|\eta_\eps''(r)|^2r |\log r|^{2\alpha}\right)dr\\
&\lesssim  \int_{\eps^{2-\mu}}^{1/2}	  r^{-3}|\log r|^{2\alpha}dr+  \int_{\eps^{2-\mu}}^{\eps } r^{-1}|\log r|^{2\alpha}\frac{dr}{r^2|\log \eps|^2}+  \int_{\eps^{2-\mu}}^{\eps } r |\log r|^{2\alpha}\frac{dr}{r^4|\log \eps|^2}\\
&\lesssim \frac{|\log \eps|^{2\alpha}}{\eps^{4-2\mu}}+\frac{|\log \eps|^{2\alpha-2}}{\eps^{4-2\mu}}\lesssim \frac{|\log \eps|^{2\alpha}}{\eps^{4-2\mu}}.
\end{align*}
The proof of the Lemma is now complete. 
\end{proof}

We now prove the desired upper bound in Theorem \ref{main}. For simplicity, we treat the case  {$N=3$}. The general case $N\ge 4$ can be handled as in \cite{ignat2025} by using Lemma \ref{lower.bound.triangle}. 
 
We consider $\tilde u_2$ as in \eqref{u2alpha}, with $\alpha>1/2,$
and the corresponding $u_\eps$ constructed previously 
\[
u_\eps(x)= \tilde u_2(|x|) \eta_\eps(|x|)\psi(|x|).
\]
 Let us denote by $\Pi_hu_\eps$ the $H_0^1$-projection of $u_\eps$ onto the space $V_h$. In view of the classical estimate \eqref{fem.error.sobolev} it satisfies
\[
\|\nabla (\Pi_h u_\eps -u_\eps)\|_{L^2(\Omega)}\lesssim h \|u_\eps \|_{H^2(\Omega)}. 
\]
Since $\Pi_h u_\eps$ is the projection of $u_\eps$ on $V_h$ it follows that
\[
\int_{\Omega}|\nabla (\Pi_h u_\eps)|^2dx=\int _\Omega |\nabla u_\eps |^2dx-\int_{\Omega}|\nabla (u_\eps -\Pi_h u_\eps)|^2dx.
\]

Using that for any $\beta>0$, $y^2\geq \frac{x^2}{1+\beta}-\frac 1\beta |x-y|^2$ and the Hardy inequality we obtain
\begin{align*}
\int_{\Omega}\frac{|\Pi_h u_\eps|^2}{|x|^2}dx&\geq \frac 1{1+\beta}\int_{\Omega}\frac{u_\eps^2}{|x|^2}dx-\frac 1\beta \int_{\Omega}\frac{|u_\eps-\Pi_hu_\eps|^2}{|x|^2}dx\\
&\geq  \frac 1{1+\beta}\int_{\Omega}\frac{u_\eps^2}{|x|^2}dx-\frac{1}{\beta\Lambda_N}\int_{\Omega}|\nabla (u_\eps -\Pi_h u_\eps)|^2dx.
\end{align*}
Let us denote
\[
R_{\eps,h}=\frac{\int_{\Omega}|\nabla (u_\eps -\Pi_h u_\eps)|^2dx}{\int_{\Omega} u_\eps^2 |x|^{-2}dx},\quad Q_\eps=\frac{\int _\Omega |\nabla u_\eps |^2dx}{\int_{\Omega} u_\eps^2 |x|^{-2}dx}.
\]
We know from Lemma \ref{est.u.eps} that \
\[
Q_\eps =\frac{A_\eps}{B_\eps}=\Lambda_N+O(|\log \eps|^{-2})
\] and
\[
R_{\eps,h}\lesssim \frac{h^2 \|u_\eps\|^2_{H^2(\Omega)}}{|\log \eps|^{2\alpha+1}}\lesssim \frac{h^2  }{|\log \eps|^{2\alpha+1}} \frac{|\log \eps|^{2\alpha}}{\eps^{4-2\mu}}=
\frac{h^2}{\eps^{4-2\mu} |\log \eps|}.
\]
Under the assumption that $\beta$ and $R_{\eps,h}$ are small enough we get
\begin{align*}
\Lambda_{h}(\Omega)& \leq \frac{\int_{\Omega}|\nabla (\Pi_h u_\eps)|^2dx}{\int_{\Omega}{|\Pi_h u_\eps|^2}|x|^{-2}dx}\leq 
\frac{\int _\Omega |\nabla u_\eps |^2dx}{\frac 1{1+\beta}\int_{\Omega}u_\eps^2|x|^{-2}dx-\frac1{\beta \Lambda_N}\int_{\Omega}|\nabla (u_\eps -\Pi_h u_\eps)|^2dx}   \\
&=\frac{Q_\eps }{\frac 1{1+\beta}-\frac 1{\beta\Lambda_N} R_{\eps,h}}\leq
\frac{\Lambda_N+O(|\log \eps|^{-2})}{1-\frac{\beta}{1+\beta}-\frac{1}{\beta\Lambda_N}\frac{h^2}{\eps^{4-2\mu} |\log \eps|}}.
\end{align*}
Taking 
\begin{equation}\eps^2\simeq h, \beta =h^{\mu/2},\quad \mu\in (0,1), \end{equation} 
we have that \begin{equation} R_{\eps, h}=O\Big (\frac{h^
\mu}{|\log h|}\Big ),  \end{equation}
and we get the desired estimate
\[
\Lambda_{h} \leq \Lambda_N+O(\frac 1{|\log h|^2}).
\]

\section{Proof of Theorem \ref{main.2}}
Classical arguments on the finite element approximation of eigenvalues show that
\begin{equation}
\mu_{1h}-\mu_1(\Omega) \sim \inf _{v_h\in V_h} \| \phi_1 -v_h\|_{\mathcal{H}}^2
\end{equation}
where $\phi_1$ is the corresponding first eigenfunction associated to $\mu_1(\Omega)$. In the case where $\Omega$ is the unit ball we have
\begin{equation} \phi_1(x)=\phi_1(|x|)=|x|^{-(N-2)/2}J_0(z_{0,1}|x|).
\end{equation} 
When $\Omega$  is a smooth domain containing the origin, one can show, using a cut-off argument, that $\phi_1$ exhibits the same type of singularity near the origin as in the case of the unit ball.

It remains to quantify the above distance in the particular case of the unit ball. 

For the upper bound, a computation similar to that  in Lemma \ref{est.u.eps} with $\alpha=0$ gives us \begin{align*}
  \inf _{v_h\in V_h} \| \phi_1 -v_h\|^2_\mathcal{H}&\lesssim  
  \|\phi_1-\phi_{1,\eps}\|^2_{\mathcal{H}} +\|\phi_{1,\eps }- I^h \phi_{1,\eps }\|^2_{\mathcal{H}}\\
  &\lesssim \frac {1}{|\log \eps|}+ \|\nabla (\phi_{1,\eps }- I^h \phi_{1,\eps })\|_{L^2(\Omega)}^2 \lesssim \frac {1}{|\log \eps|}+ \frac{h^2}{\eps^{4-2\mu}}.
\end{align*}
Taking $\eps^2\simeq h$, we obtain the desired upper bound. 

For the lower bound, we must show that for any  $v_h\in V_h$,
\[
 \| \phi_1 -v_h\|^2_{\mathcal{H}}\gtrsim \frac 1{|\log h|}.
\]
Using the Hardy inequality with logarithmic correction in Lemma \ref{hardy.log.2}, we obtain 
\begin{align*}
\label{}
   \| \phi_1 -v_h\|^2_{\mathcal{H}}\gtrsim \int_{\Omega}\frac{|\nabla \phi_1-\nabla v_h|^2dx}{|\log(1/|x|)|^2}=\sum_{T\in \mathcal{T}_h} \int_T \frac{|\nabla \phi_1-\nabla v_h|^2dx}{|\log(1/|x|)|^2}.
\end{align*}
It is interesting to notice that for all the simplices situated outside the ball of radius $h$, the contribution is of order $1/{|\log h|^2}$. The main contribution comes from the simplex $T_0$  containing the origin or the finitely many simplices for which the origin lies on their boundaries.

 To simplify the presentation, let us consider the case in which the origin is the incenter of $T_0$,
and a ball of radius $\rho_{T_0}\simeq h$ centered at $x=0$ is contained in $T_0$. Otherwise, we can perform a similar computation by choosing a conical subset. 

In this case, the main contribution is given by
\[
I_h= \int_{|x|\leq h} \frac{|\nabla \phi_1-\nabla v_h|^2dx}{|\log(1/|x|)|^2}\geq \inf _{A\in \rr^N} \int_{|x|\leq h} \frac{|\nabla \phi_1-A|^2dx}{|\log(1/|x|)|^2}.
\] 
Since 
\[
\nabla \phi_1(x)=\phi_1'(|x|)\frac x{|x|},
\]
 we obtain that
 \[
 |\nabla \phi_1-A|\geq ||\phi'_1(r)|-|A||. 
 \]
 We have to prove that for any $A\in \mathbb{R}^N$,
 \[
 \int_{B_h(0)} \frac{||\phi'_1(|x|)|- |A||^2}{|\log(1/|x|)|^2}dx \gtrsim \frac 1{|\log h|}.
 \]
 After integrating in the angular variables, this reduces to showing that for any constant $C_A\ge 0$,
 \[
 g(C_A)=\int _0^h \frac {r^{N-1}(|\phi_1'(r)|-C_A)^2}{|\log(r)|^2}dr \gtrsim \frac 1{|\log h|},
 \]
 where $\phi_1(r)=r^{-(N-2)/2}J_0(z_{0,1}r)$. 

 Denoting $g(C_A)=C_A^2 a_h -2C_A b_h+c_h$, we know that the minimum is attained at $C_A=b_h/a_h$, and 
 \[
 \min _{C_A\ge 0} g(C_A) = c_h-\frac{b_h^2}{a_h}.
  \]
 We use the following expansion as $h\rightarrow 0$, valid for $\alpha\ge -1$:
 \[
 \int_0^h \frac{r^\alpha}{|\log r|^2}dr\simeq 
 \begin{cases}
 	\frac{h^{\alpha+1}}{|\log h|^2},& \alpha>-1,\\
 		\frac{1}{|\log h|},& \alpha=-1.
 \end{cases} 
 \]
 Since $\phi_1'(r)\simeq r^{-N/2}$ as $r\rightarrow 0$, we have $c_h = \int_0^h \frac{r^{N-1}(\phi_1'(r))^2}{|\log r|^2}dr \simeq \int_0^h \frac{r^{-1}}{|\log r|^2} dr \simeq \frac{1}{|\log h|}$. Similarly, we obtain $a_h \simeq \frac{h^N}{|\log h|^2}$ and $b_h \simeq \frac{h^{N/2}}{|\log h|^2}$.
 Then 
 \[
  \min _{C_A\ge 0 } g(C_A) \simeq \frac 1{|\log h|} - \frac{(h^{N/2}/|\log h|^2)^2}{h^N/|\log h|^2} = \frac 1{|\log h|} - \frac{1}{|\log h|^2} \simeq \frac{1}{|\log h|}.
 \]
  
 \section{Proof of Theorem \ref{main.3}} We consider the case where the domain $\Omega$ is the unit ball, since the singularity analysis near the origin is the same for all domains containing $x=0$.
 
In the present case,  the norm introduced by the bilinear form is equivalent to the $H^1_0(\Omega)$-one.  Thus, it is sufficient to provide sharp upper bounds on \[
 d(\phi_1, V_h)\simeq \inf_{v_h\in V_h} \|\nabla \phi_1 -\nabla v_h \|_{H^1_0(B_1)},
 \]
 where $\phi_1$ is the first eigenfunction.

  Let us set \[
 m=\sqrt{\Lambda_N-\Lambda}\in \left(0,\frac{N-2}2\right].\]  Recall  that  
  in the case of the ball, the first eigenfunction can be explicitly computed  (see \cite{MR1760280}), $\phi_1(x)=\phi_1(|x|)$,
 \[
 \phi_1(r)=r^{-N/2+1} J_{m}(j_{m,1}r)\simeq r^{-N/2+1}\left(\frac{1}{\Gamma(m+1)} \left( \frac{r}{2} \right)^m - \frac{1}{\Gamma(m+2)} \left( \frac{r}{2} \right)^{m+2}\right), r\rightarrow 0,
 \]
 where $J_m$ is the Bessel function of order $m$ and $j_{m,1}$ is its  first zero.  
 
 This eigenfunction $\phi_1$ belongs to $H^{1+s}(\Omega)$ for all $0<s<m$.  If $m=(N-2)/2$ then  $\phi_1\in C^\infty(\Omega)$. 

Let us prove the upper bound.  Let us assume without loss of generality that there exists a simplex $T_0$ that contains the origin and all its vertices have the same distance to the origin. 
We take a function $v_h$ that is constant in the simplex $T_0$ and interpolate $\phi_1$ at the other nodes.  We quantify the interpolation error by splitting the domain into a small ball around the origin, where the singular behavior dominates, and its complement, where standard $H^2$-based estimates apply. We split the integral into two parts:
\[
\|\nabla \phi_1 -\nabla v_h \|_{L^2(\Omega)}^2\leq \int_{|x|\lesssim h}|\nabla \phi_1|^2dx +h^2 \int_{|x|\gtrsim h} |D^2\phi_1|^2dx:=I_1+I_2.
\]
When $m<\frac N2-1$ we have 
\[
I_1\leq \int_0^h r^{N-1}r^{2(-N/2+m)}dr\simeq h^{2m}
\] 
while for $m=\frac N2 -1$ we have $I_1\simeq h^N$ since $
\phi_1$ is $C^\infty$ in this case.  

For the second term, we first consider the case $N\geq 5$. When $m=\frac{N-2}2$ clearly we have $I_2\lesssim h^2$.  For $0<m<\frac{N-2}2$ a similar argument as for $I_1$ yields
\[
I_2\lesssim h^2\int_h^1 r^{N-1}r^{2(-N/2-1+m)}dr=h^2\int_h^1 r^{2m-3}dr  \lesssim   \begin{cases}
h^{2m},& 0<m<1,\\
h^2|\log h|, & m=1,\\
h^2, &m\in (1, \frac{N-2}2]. 
  \end{cases}
\]
Putting together the results for $I_1$ and $I_2$ leads to the desired result.

When $N=4$ we obtain, in a similar way, that $I_2\leq h^2$ for $m=\frac{N-2}2=1$ and \[
I_2\leq h^2\int_h^1 r^{2m-3}dr\simeq h^{2m}, \ 0<m<1.
\]

In dimension $N=3$ we get $I_2\leq h^2$ when $m=\frac 12=\frac{N-2}2$ and 
 \[
 I_2\leq h^2 \int_{|x|\gtrsim h} |D^2\phi_1|^2dx
\leq h^2\int_h^1 r^{2m-3}dr \simeq h^{2m}, \ 0<m<1/2.
\]

To obtain the lower bound in \eqref{est.hardy.3} we divide the integral into two parts, the first one as in the proof of Theorem \ref{main.2} and 
the second one, integrating outside the ball of radius $h$:
\[
 \|\nabla \phi_1 -\nabla v_h \|_{L^2(\Omega)}^2=\int_{|x|<h}|\nabla \phi_1 -\nabla v_h|^2dx+\int_{h<|x|<1}|\nabla \phi_1 -\nabla v_h|^2dx=I_1+I_2.
\]
Let us consider the case $N\geq 5$ since the others are similar.  For $I_1$ the same arguments as in the proof of Theorem \ref{main.2}
 give us that $I_1\gtrsim h^{2m}$ for $0<m<\frac{N-2}2$ or $I_1\gtrsim h^{2N}$ if $m=\frac{N-2}2$. 
Using Lemma \ref{lower.bound.triangle} as in \cite{ignat2025} we get 
 \[
 I_2\gtrsim  \begin{cases}
h^{2m},& 0<m<1,\\
h^2|\log h|, & m=1,\\
h^2, &m>1. 
  \end{cases}
 \]
 The proof is now complete.

\section{Conclusions and Open Problems}

In this paper, we have analyzed the finite element approximation of the best Hardy constant in bounded domains containing the origin for dimensions $N \geq 3$. Despite the absence of minimizers for the Hardy inequality in the standard Sobolev space $H^1_0(\Omega)$, we have rigorously established that the first eigenvalue of the corresponding discrete eigenvalue problem converges to the continuous Hardy constant as the mesh size $h \to 0$.

Our main result provides an explicit convergence rate of order $1/|\log h|^2$, independent of the spatial dimension, and reflects the singular nature of the underlying functional inequality. The analysis demonstrates how the finite element method is capable of capturing the concentration phenomena inherent in the minimization sequences that saturate the Hardy inequality.

We have also analyzed similar approximation issues for two closely related eigenvalue problems.

Several open problems naturally arise from this study:
\begin{itemize}
\item{\it The $2$-dimensional case:} As mentioned in the introduction, the two-dimensional case is critical and requires an inverse-square logarithmic correction to the Hardy inequality, as shown in \cite{MR1938711}. This case would therefore require a more detailed analysis.

\item {\it Other Hardy inequalities:} The literature on Hardy inequalities is extensive; we refer the reader to the monograph \cite{peral} for a comprehensive overview. The finite element approximation techniques developed in this paper could, in principle, be systematically applied to the many existing variants of the Hardy inequality. 
\item {\it Weighted spectral problems:} Using the fact  that 
 \[  \mu_1(\Omega)=\min _{u\in \mathcal{H}}\frac{\|u\|^2_{ \mathcal{H}} }{\|u\|^2_{L^2(\Omega)}}=
\min _{v\in W^{1,2}_{0}(|x|^{-(N-2)}dx,\Omega)} \frac{\int _\Omega  |x|^{-(N-2)} |\nabla v|^2dx}{\int_ \Omega |x|^{-(N-2)} v^2dx},
\]
one can consider introducing the finite element approximation in this weighted setting 
 \[\tilde \mu_{1h}(\Omega)=
\min _{v\in V_h} \frac{\int _\Omega  |x|^{-(N-2)} |\nabla v|^2dx}{\int_ \Omega |x|^{-(N-2)} v^2dx}
\]
and analyze the approximation rates.

	\item {\it Higher-Order Finite Elements:}
It would be of interest to investigate whether similar convergence rates can be established for higher-order finite element spaces and other finite-element variants, like discontinuous Galerkin approximations.
	\item	 {\it  Adaptive Mesh Refinement:}
Given the concentration of minimizing sequences near the singularity at the origin, adaptive mesh refinement strategies could improve the numerical approximation. Theoretical analysis of such adaptive methods in this singular context remains an open challenge.
	\item	 {\it Extension to Non-Radial Settings and General Domains:}
While our analysis is presented for general bounded domains containing the origin, extensions to more complex geometries or domains with additional singularities could reveal new phenomena. The same can be said in the context of multi-polar Hardy inequalities (see \cite{cazacu2013improved} for the analysis in the continuous setting).
	\item	{\it Other Singular Inequalities:}
A natural direction for future work is to study the finite element approximation of best constants in other critical inequalities, such as Rellich or Hardy-Sobolev inequalities, where similar non-attainability issues arise. Similar issues also arise in the fractional setting, see \cite{de2024fractional}.

\end{itemize}
We hope that this contribution stimulates further research on the numerical analysis of critical inequalities and the development of numerical methods tailored for singular variational problems.

 \subsection*{Acknowledgements}
L. I.  Ignat  was supported by a grant of the Ministry of Research, Innovation, and Digitization, CCCDI -
UEFISCDI, project number ROSUA-2024-0001, within PNCDI IV.

E. Zuazua was funded by the European Research Council (ERC) under the European Union's Horizon 2030 research and innovation programme (Grant No. 101096251-CoDeFeL), the Alexander von Humboldt-Professorship program, the ModConFlex Marie Curie Action, HORIZON-MSCA-2021-dN-01, the Transregio 154 Project of the DFG, grants PID2020-112617GB-C22 and TED2021131390B-I00 of the AEI (Spain),  AFOSR Proposal 24IOE027, and Madrid Government-UAM Agreement for the Excellence of the University Research Staff in the context of the V PRICIT (Regional Programme of Research and Technological Innovation).

\section{Appendix}

\begin{lemma}\label{hardy.log.2}
Let $\Omega$ be a  smooth domain such that 
$0\in \Omega \subset \overline{\Omega}\subset B_R$ for some $R>0$. There exists a positive constant $C=C(N,\Omega)$ such that for all $u\in \mathcal{H}$ we have
\begin{equation}
\label{h.log}
\|u\|_{\mathcal{H}}^2 \geq C(N,\Omega)\int_{\Omega}\frac{|\nabla u|^2dx}{|\log(R/|x|)|^2}.
\end{equation}
\begin{proof}
	The inequality holds for $u\in C_c^\infty(\Omega)$ as proved in  \cite[Th.~2.5.2, p.~25]{peral}
\[
\|u\|_{\mathcal{H}}^2 = \int_\Omega |\nabla u|^2 dx- \Lambda_N\int_{\Omega}\frac{u^2}{|x|^2}dx\geq C \int_{\Omega} |\nabla u|^2 \left(\log\frac R {|x|} \right)^{-2}dx.
\]
	 By density, the inequality extends to functions $u$ in $\mathcal{H}$. 
\end{proof}

\end{lemma}

\begin{lemma}
\label{logth}
For any $T\in \th$ the following holds:
\begin{equation}
\label{est.log.below}
  \int_T \frac{dx}{|\log |x||^2}\gtrsim \frac{|T|}{|\log h|^2}.
\end{equation}
\end{lemma}

The proof is elementary and is left to the readers.
 \bibliographystyle{abbrv}

\begin{thebibliography}{10}

\bibitem{MR1938711}
Adimurthi and K.~Sandeep.
\newblock Existence and non-existence of the first eigenvalue of the perturbed
  {H}ardy-{S}obolev operator.
\newblock {\em Proc. Roy. Soc. Edinburgh Sect. A}, 132(5):1021--1043, 2002.

\bibitem{MR2754218}
P.~F. Antonietti and A.~Pratelli.
\newblock Finite element approximation of the {S}obolev constant.
\newblock {\em Numer. Math.}, 117(1):37--64, 2011.

\bibitem{buttazzo}
H.~Attouch, G.~Buttazzo and G.~Michaille.
\newblock {\em Variational analysis in {S}obolev and {BV} spaces}, volume~6 of
  {\em MPS/SIAM Series on Optimization}.
\newblock Society for Industrial and Applied Mathematics (SIAM), Philadelphia,
  PA; Mathematical Programming Society (MPS), Philadelphia, PA, 2006.
\newblock Applications to PDEs and optimization.

\bibitem{MR1115240}
I.~Babu\v{s}ka and J.~Osborn.
\newblock Eigenvalue problems.
\newblock In {\em Handbook of numerical analysis, {V}ol. {II}}, volume~II of
  {\em Handb. Numer. Anal.}, pages 641--787. North-Holland, Amsterdam, 1991.

\bibitem{MR2373954}
S.~C. Brenner and L.~R. Scott.
\newblock {\em The mathematical theory of finite element methods}, volume~15 of
  {\em Texts in Applied Mathematics}.
\newblock Springer, New York, third edition, 2008.

\bibitem{MR3528527}
C.~Cazacu and D.~Krej{\v c}i{\v r}{\'i}k.
\newblock The {H}ardy inequality and the heat equation with magnetic field in
  any dimension.
\newblock {\em Comm. Partial Differential Equations}, 41(7):1056--1088, 2016.

\bibitem{cazacu2013improved}
C.~Cazacu and E.~Zuazua.
\newblock Improved multipolar hardy inequalities.
\newblock In {\em Studies in phase space analysis with applications to PDEs},
  pages 35--52. Springer, 2013.

\bibitem{MR3405533}
F.~Chatelin.
\newblock {\em Spectral approximation of linear operators}, volume~65 of {\em
  Classics in Applied Mathematics}.
\newblock Society for Industrial and Applied Mathematics (SIAM), Philadelphia,
  PA, 2011.

\bibitem{MR2078192}
E.~Colorado and I.~Peral.
\newblock Eigenvalues and bifurcation for elliptic equations with mixed
  {D}irichlet-{N}eumann boundary conditions related to
  {C}affarelli-{K}ohn-{N}irenberg inequalities.
\newblock {\em Topol. Methods Nonlinear Anal.}, 23(2):239--273, 2004.

\bibitem{de2024fractional}
N.~De~Nitti and S.~M. Djitte.
\newblock Fractional {Hardy}--{Rellich} inequalities via integration by parts.
\newblock {\em Nonlinear Analysis}, 243:113478, 2024.

\bibitem{della2023finite}
F.~Della~Pietra, G.~Fantuzzi, L.~I. Ignat, A.~L. Masiello, G.~Paoli and
  E.~Zuazua.
\newblock Finite element approximation of the {H}ardy constant.
\newblock {\em J. Convex Anal.}, 31(2):497--523, 2024.

\bibitem{MR3823299}
D.~E. Edmunds and W.~D. Evans.
\newblock {\em Spectral theory and differential operators}.
\newblock Oxford Mathematical Monographs. Oxford University Press, Oxford,
  second edition, 2018.

\bibitem{ignat2025}
L.~I. Ignat and E.~Zuazua.
\newblock Optimal convergence rates for the finite element approximation of the
  {Sobolev} constant.
\newblock {\em arXiv:2504.09637}, 2025.

\bibitem{peral}
I.~Peral~Alonso and F.~Soria~de Diego.
\newblock {\em Elliptic and parabolic equations involving the {Hardy}-{Leray}
  potential}, volume~38 of {\em De Gruyter Ser. Nonlinear Anal. Appl.}
\newblock Berlin: De Gruyter, 2021.

\bibitem{raviart}
P.-A. Raviart and J.-M. Thomas.
\newblock {\em Introduction \`a l'analyse num\'{e}rique des \'{e}quations aux
  d\'{e}riv\'{e}es partielles}.
\newblock Collection Math\'{e}matiques Appliqu\'{e}es pour la Ma\^{i}trise.
  [Collection of Applied Mathematics for the Master's Degree]. Masson, Paris,
  1983.

\bibitem{MR2960851}
J.~L. V\'azquez and N.~B. Zographopoulos.
\newblock Functional aspects of the {H}ardy inequality: appearance of a hidden
  energy.
\newblock {\em J. Evol. Equ.}, 12(3):713--739, 2012.

\bibitem{MR1760280}
J.~L. Vazquez and E.~Zuazua.
\newblock The {H}ardy inequality and the asymptotic behaviour of the heat
  equation with an inverse-square potential.
\newblock {\em J. Funct. Anal.}, 173(1):103--153, 2000.

\end{thebibliography}

\end{document}